\documentclass[11pt,a4paper]{article}%
\usepackage{amsmath, amsfonts, color, amsthm, latexsym}
\usepackage{amsmath}
\usepackage{amsfonts}
\usepackage{amssymb}
\usepackage{graphicx}%
\providecommand{\U}[1]{\protect\rule{.1in}{.1in}}
\newtheorem{theorem}{Theorem}[section]
\newtheorem{proposition}[theorem]{Proposition}
\newtheorem{corollary}[theorem]{Corollary}
\newtheorem{example}[theorem]{Example}
\newtheorem{remark}[theorem]{Remark}
\newtheorem{remarks}[theorem]{Remarks}

\begin{document}

\title{\textsc{Lineability of summing sets of homogeneous polynomials}}
\author{Geraldo Botelho\thanks{Corresponding author. Partially supported by CNPq Grant 202162-1006-0.}\,\,,
M\'ario C. Matos and Daniel Pellegrino\thanks{Partially supported by CNPq Grants 471054/2006-2 and
308084/2006-3.\hfill\newline2000 Mathematics Subject Classification. Primary 46G25; Secondary
47B10.\hfill\newline Keywords: Homogeneous polynomials, absolute summability, lineability, spaceability.} }
\date{}
\maketitle

\begin{abstract}
Given a continuous $n$-homogeneous polynomial $P\colon E\longrightarrow F$
between Banach spaces and $1\leq q\leq p<\infty$, in this paper we investigate
some properties concerning lineability and spaceability of the $(p;q)$-summing
set of $P$, defined by $S_{p;q}(P)=\{a\in E:P\mathrm{~is~}%
(p;q)\mathrm{-summing~at~}a\}$.

\end{abstract}

\vspace*{-1.0em}

\section*{Introduction}

Several different generalizations of absolutely summing linear operators to
homogeneous polynomials have already been studied and besides its intrinsic
mathematical interest, this line of investigation has interesting applications
(for example, in the study of convolution equations: a nonlinear concept
related to absolutely summing operators is used in \cite{Favaro, MM} to
generalize results from \cite{Mal, Mar, MCM}). One of the possible polynomial
extensions of the concept of absolutely summing operator is the class of
homogeneous polynomials which are absolutely summing at a given point of the
domain, which was introduced by M. C. Matos \cite{Nach} and developed in
\cite{BBDP, BBJP, Erhan, PArchiv, irish}. This class is interesting, among
other reasons, because it is suitable to extend the theory to arbitrary
nonlinear mappings in the following fashion: given $1\leq q\leq p<\infty$ and
a mapping $f\colon E\longrightarrow F$ between Banach spaces, we say that $f$
is $(p;q)$-summing at a point $a\in E$ if $(f(a+x_{j})-f(a))_{j=1}^{\infty}$
is absolutely $p$-summable in $F$ whenever $(x_{j})_{j=1}^{\infty}$ is
unconditionally $q$-summable in $E$. \newline\indent The $(p;q)$-summing set
of the mapping $f$ is defined by
\[
S_{p;q}(f)=\{a\in E:f\mathrm{~is~}(p;q)\mathrm{-summing~at~}a\}.
\]
If $p=q$ we just say that $f$ is $p$-summing at $a$ and simply write
$S_{p}(f)$. Letting $E$ be infinite-dimensional, the following questions are
natural: Is $S_{p;q}(f)$ non-empty? If yes, is it a linear subspace of $E$? If
yes, is it $\neq\{0\}$? If yes, is it infinite-dimensional? Closed? The whole
space? If no, does it contain a linear subspace $G$ of $E$? If yes, does it
contain an infinite-dimensional subspace? Closed?\newline\indent In this paper
we address the above questions for homogeneous polynomials. In Section 3 we
characterize the summing set of a polynomial in terms of the summability of
its derivatives at the origin and we study the cases of polynomials $P\colon
E\longrightarrow F$ for which $S_{p;q}(P)=\emptyset$, $S_{p;q}(P)=\{0\}$ and
$S_{p;q}(P)=E$. In section 4 we identify the precise point where linearity is
lost. More specifically, we prove that the summing set of either a
2-homogeneous polynomial or a scalar-valued 3-homogeneous polynomial is a
linear subspace; whereas the summing set of either a vector-valued
3-homogeneous polynomial or a scalar-valued $n$-homogeneous polynomials,
$n\geq4$, may fail to be a linear subspace. We also prove that, even when the
summing set is an infinite-dimensional subspace, it may fail to contain a
closed infinite-dimensional subspace. In a final result we prove that,
regardless of the positive integer $n\geq2$, any finite-dimensional subspace
of $L_{2}$ is the summing set of a certain $n$-homogeneous polynomial.

It is worth mentioning that this line of investigation, for different sorts of sets, has been previously
explored. For example, given a continuous homogeneous polynomial $P \colon E \longrightarrow F$, the sets
$$C(P)=\{a\in E; P\text{ is weakly sequentially continuous at }a \}{\rm ~and}$$
$$
c_{w}(P)=\{a\in E;P\text{ is weakly continuous on bounded sets at }a\} $$ were investigated in
\cite{Aron-Dimant, neusa} and \cite{Boyd-Ryan}, respectively (see Section 3).


\section{Background and notation}

Throughout this paper $E$ and $F$ will stand for Banach spaces over
$\mathbb{K}=\mathbb{R}$ or $\mathbb{C}$, $E^{\prime}$ is the dual of $E$ and
$n$ will always be a positive integer. By $\mathcal{P}(^{n}E;F)$ we denote the
Banach space of all continuous $n$-homogeneous polynomials from $E$ to $F$
with the usual sup norm. If $F=\mathbb{K}$ we simply write $\mathcal{P}%
(^{n}E)$. Given $P\in\mathcal{P}(^{n}E;F)$, $a\in E$ and $k\in\{1,\ldots,n\}$,
denoting by $\check{P}$ the symmetric $n$-linear mapping associated to $P$, as
usual we define:
\[
\hat{d}^{k}P(a)\colon E\longrightarrow F~:~\hat{d}^{k}P(a)(x)=\frac
{n!}{(n-k)!}\check{P}(a^{n-k},x^{k})~;~d^{k}P(a)\in\mathcal{P}(^{k}E;F);
\]%
\[
\hat{d}^{k}P\colon E\longrightarrow\mathcal{P}(^{k}E;F)~:~a\longrightarrow
d^{k}P(a)~;~\hat{d}^{k}P\in\mathcal{P}(^{n-k}E;\mathcal{P}(^{k}E;F)).
\]
Here, $(a^{n-k},x^{k})$ means $(a,\ldots,a,x,\ldots,x)$, where $a$ appears
$(n-k)$ times and $x$ appears $k$ times. For $k=1$ we write $\hat{d}P(a)$ and
$\hat{d}P$ instead of $\hat{d}^{1}P(a)$ and $\hat{d}^{1}P$. For the general
theory of multilinear mappings, homogeneous polynomials and holomorphic
mappings we refer to \cite{din, mu}.\newline\indent Let $p\geq1$. By $\ell
_{p}(E)$ we mean the Banach space of all absolutely $p$-summable sequences
$(x_{j})_{j=1}^{\infty}$, $x_{j}\in E$ for all $j$, with the norm $\Vert
(x_{j})_{j=1}^{\infty}\Vert_{p}=\left(  \sum_{j=1}^{\infty}\Vert x_{j}%
\Vert^{p}\right)  ^{1/p}$. $\ell_{p}^{w}(E)$ denotes the Banach space of all
sequences $(x_{j})_{j=1}^{\infty}$, $x_{j}\in E$ for all $j$, such that
$(\varphi(x_{j}))_{j=1}^{\infty}\in\ell_{p}$ for every $\varphi\in E^{\prime}$
with the norm
\[
\Vert(x_{j})_{j=1}^{\infty}\Vert_{w,p}=\sup\{\Vert(\varphi(x_{j}%
))_{j=1}^{\infty}\Vert_{p}:\varphi\in E^{\prime},\Vert\varphi\Vert\leq1\}.
\]
$\ell_{p}^{u}(E)$ is the closed subspace of $\ell_{p}^{w}(E)$ formed by the
sequences $(x_{j})_{j=1}^{\infty}$ satisfying $\lim_{k\rightarrow\infty}%
\Vert(x_{j})_{j=k}^{\infty}\Vert_{w,p}=0$. Such sequences are called
unconditionally $p$-summable.\newline\indent According to the definition
stated in the introduction, given $1\leq q\leq p<\infty$, a mapping $f\colon
E\longrightarrow F$ is $(p;q)$-summing at $a\in E$ if $(f(a+x_{j}%
)-f(a))_{j=1}^{\infty}\in\ell_{p}(F)$ whenever $(x_{j})_{j=1}^{\infty}\in
\ell_{q}^{u}(E)$. For a polynomial $P\in\mathcal{P}(^{n}E;F)$, by
\cite[Proposition 2.4]{ma96} we know that $P$ is $(p;q)$-summing at the origin
if and only if $(P(x_{j}))_{j=1}^{\infty}\in\ell_{p}(F)$ whenever
$(x_{j})_{j=1}^{\infty}\in\ell_{q}^{w}(E)$. Making $n=1$ we recover the
classical ideal of $(p;q)$-summing linear operators. The space of all
$(p;q)$-summing linear operators from $E$ to $F$ will be denoted by $\Pi
_{p;q}(E;F)$. For the theory of absolutely summing linear operators we refer
to \cite{djt}.\newline\indent According to \cite{ags, gq} and others, a subset
$A$ of a topological vector space $E$ is said to be:\newline$\bullet$
\textit{$n$-lineable} if $A\cup\{0\}$ contains an $n$-dimensional linear
subspace of $E$.\newline$\bullet$ \textit{lineable} if $A\cup\{0\}$ contains
an infinite-dimensional linear subspace of $E$.\newline$\bullet$
\textit{spaceable} if $A\cup\{0\}$ contains a closed infinite-dimensional
linear subspace of $E$.

\section{Known related facts}

Several known results and examples can be rephrased in the language of summing
sets. We list some of them in this section in order to inform the reader about
the state of the art and for further reference.

\begin{example}
\label{Example 2.1}\rm\cite[Example 3.2 and Theorem 6.3]{Nach} Let $E$ be an infinite-dimensional Banach
space, $n \geq2$, $p \geq1$ and $\varphi\in E^{\prime}$, $\varphi\neq0$. If $P \colon E \longrightarrow E$ is
the $n$-homogeneous polynomial defined by $P(x) = \varphi(x)^{n-1}x$, then $S_{p}(P) = ker(\varphi)$.
\end{example}

\begin{proposition}
\label{Proposition 2.2}{\rm\cite[Corollary 3.6]{BBDP}} Let $P\in \mathcal{P}(^{n}E;F)$, $n\geq2$. If $P$ is
$(p;q)$-summing at $a\in E$, then $P$ is $(p;q)$-summing at $\lambda a$ for every $\lambda\in\mathbb{K}$. In
consequence, $S_{p;q}(P)$ is either empty, or $\{0\}$ or 1-lineable. In particular $S_{p;q}(P)\neq\emptyset$
if and only if $0\in S_{p;q}(P)$.
\end{proposition}

\begin{proposition}
\label{Proposition 2.3}{\rm\cite[Theorem 3.10]{am}} Let $n\geq2$ and $E$ be a Banach space. For every
$P\in\mathcal{P}(^{n}E)$, $0\in S_{1}(P)$.
\end{proposition}

Remember that a Banach space $E$ has the \textit{Orlicz property} if the
identity operator on $E$ is $(2;1)$-summing. Spaces of cotype 2 have the
Orlicz property.

\begin{proposition}
\label{orlicz}{\rm\cite[Proposition 2.9]{ma96}} Let $n \geq2$, $E$ be a Banach space with the Orlicz
property, $F$ be an arbitrary Banach space and $P \in\mathcal{P}(^{n}E;F)$. Then $0 \in S_{1}(P)$.
\end{proposition}

If a mapping $f \colon E \longrightarrow F$ is such that $S_{p;q}(f) = E$, $f$
is called \textit{everywhere $(p;q)$-summing}. Several such cases are known,
for example:

\begin{proposition}
~\newline{\rm(a) \cite[Corollary 4.2]{Nach}} Let $E$ be an $\mathcal{L}_{1}$-space and $F$ be an
$\mathcal{L}_{2}$-space. Then
$S_{1}(f)=E$ for every analytic mapping $f\colon E\longrightarrow F$%
.\newline{\rm (b) \cite[Theorem 5.2]{irish}} Let $F$ be a Banach space of cotype $q$ and $E$ be an arbitrary
Banach space. Then $S_{q;1}(f)=E$ for every analytic mapping $f\colon E\longrightarrow F$.
\newline{\rm(c)
\cite[Theorem 2.4]{irish}} Let $n\in\mathbb{N}$. A Banach space $E$ has cotype
$q>2$ if and only if $S_{q;1}(P)=E$ for every $P\in\mathcal{P}(^{n}%
E;E)$.
\newline{\rm(d) \cite[Theorem 2.10]{irish}} If $E$ is an
$\mathcal{L}_{\infty}$-space and $F$ has cotype $q$, then $S_{q;2}(P)=E$ for
every $P\in\mathcal{P}(^{n}E;F)$.
\end{proposition}

\section{General examples and results}

Our first three examples show that, though apparently similar at first glance, the theory of summing sets
$S_{p;q}(P)$ is actually quite different from that of the sets $C(P)$ of \cite{Aron-Dimant} (for the
definition see the Introduction). Recall that a polynomial $P\in\mathcal{P}(^{n}E;F)$ is
\textit{$p$-dominated} \cite[Definition 3.2]{ma96}, $1\leq p<+\infty$, if
$(P(x_{j}))_{j=1}^{\infty}\in\ell_{\frac {p}{n}}(F)$ whenever $(x_{j})_{j=1}^{\infty}\in\ell_{p}^{w}(E)$.

\begin{example}
\rm Let $P\in\mathcal{P}(^{n}E;F)$ be a $p$-dominated polynomial. On the one hand, $S_{p}%
(P)=C(P)=E.$ Indeed, the fact that $P$ can be written as $P = Q \circ u$ where $u$ is an absolutely
$p$-summing linear operator, hence completely continuous, shows that $C(P)=E.$ The fact that $S_{p}(P)=E$ is
a straightforward consequence of \cite[Theorem 3((iii) and (v))]{Erhan}. On the other hand, in general
$S_{q}(P)\neq C(P)$ for $q<p$. For example, $S_1(v) = \emptyset \neq E =C(v)$ for any absolutely 2-summing
non-absolutely 1-summing linear operator $v$ on $E$.
\end{example}

In particular, from the previous example we conclude that when $E$ is either an $\mathcal{L}_{\infty}$-space,
the disc algebra $\mathcal{A}$ or the Hardy space $H^{\infty},$ then $S_{p}(P)=C(P)=E$ for every
$P\in\mathcal{P}(^{2}E)$ and every $p\geq2$. In fact, from \cite[Proposition 2.1]{Pams} we know that every
such $P$ is $p$-dominated ($p\geq2$).

\begin{example}
\rm The polynomial
\[
P\colon\ell_{2}\longrightarrow\mathbb{R}~;~P((\alpha_{j})_{j=1}^{\infty}%
)=\sum_{j=1}^{\infty}\alpha_{j}^{2}.
\]
is absolutely $1$-summing (and hence $0\in S_{p}(P)),$ but $C(P)=\emptyset$ because
$(a+e_{j})_{j=1}^{\infty}$ is weakly convergent to $a$ but $(P(a+e_{j}))_{j=1}^{\infty}$ fails to be norm
convergent to $P(a)$.
\end{example}

\begin{example}
\rm If $\varphi\in E^{\prime},$ we know from Example \ref{Example 2.1} that $P(x)=\varphi(x)x$ is so that $
S_{1}(P)=Ker(\varphi)$. On the other hand, if $E=\ell_{1}$ we have $C(P)=\ell_{1}$ since $\ell_{1}$ has the
Schur property.
\end{example}

All summing sets of finite type polynomials and 1-summing sets of nuclear polynomials can be easily
described:

\begin{example}
\label{tipo finito}\rm Let $P\in\mathcal{P}(^{n}E;F)$ be a polynomial of finite type, that is
$P(x)=\sum_{j=1}^{k}\varphi_{j}(x)^{n}b_{j}$, where $k\in\mathbb{N}$, $\varphi_{1},\ldots,\varphi_{k}\in
E^{\prime}$ and $b_{1},\ldots,b_{k}\in F$. An easy adaptation of the proof of \cite[Lemma 6.2]{Nach} shows
that $S_{p;q}(P)=E$ for every $1\leq q\leq p$.
\end{example}

\begin{example}
\rm A polynomial $P\in\mathcal{P}(^{n}E;F)$ is said to be \textit{nuclear} \cite[Definition 2.9]{din} if
there exist $(\lambda
_{j})_{j=1}^{\infty}\in\ell_{1}$ and bounded sequences $(\varphi_{j}%
)_{j=1}^{\infty}$ in $E^{\prime}$ and $(b_{j})_{j=1}^{\infty}$ in $F$ such
that
\[
P(x)=\sum_{j=1}^{\infty}\lambda_{j}\varphi_{j}(x)^{n}b_{j}~\mathrm{for~every~}%
x\in E.
\]
It is not difficult to show that $S_{1}(P)=E$ for every nuclear polynomial $P\in\mathcal{P}(^{n}E;F)$.
\end{example}

Let us see that the summing set may be empty even for simple nonlinear
mappings. For $p=2$ the scalar-valued case is enough:

\begin{example}
\label{Example 3.1} \rm Consider the $n$-homogeneous polynomial (n\textrm{$\geq2)$}
\[
P\colon\ell_{2}\longrightarrow\mathbb{K}~;~P((\alpha_{j})_{j=1}^{\infty}%
)=\sum_{j=1}^{\infty}\alpha_{j}^{n}.
\]
Let $2\leq q\leq p$. Since $(e_{j})_{j=1}^{\infty}\in\ell_{2}^{w}(\ell _{2})\subseteq\ell_{q}^{w}(\ell_{2})$
and $P(e_{j})=1$ for every $j$, we have that $P$ is not $(p;q)$-summing at $0$. So, $S_{p;q}(P)=\emptyset$ by
Proposition \ref{Proposition 2.2} for every $2\leq q\leq p$. In particular $S_{p}(P)=\emptyset$ for every
$p\geq2$.
\end{example}

For $p = 1$, Proposition \ref{Proposition 2.3} forces a vector-valued example:

\begin{example}
\label{Example 3.2}\rm Consider the $n$-homogeneous polynomial
\[
P \colon c_{0} \longrightarrow c_{0} ~;~ P((\alpha_{j})_{j=1}^{\infty}) =
(\alpha_{j}^{n})_{j=1}^{\infty}.
\]
Let $1 \leq q \leq p$. Since $(e_{j})_{j=1}^{\infty}\in\ell_{1}^{w}(c_{0}) \subseteq\ell_{q}^{w}(c_{0})$ and
$\|P(e_{j})\| = \|e_{j}\| = 1$ for every $j$, we have that $P$ is not $(p,q)$-summing at $0$. So, $S_{p;q}(P)
= \emptyset$ by Proposition \ref{Proposition 2.2}. for every $1 \leq q \leq p$. In particular $S_{p}(P) =
\emptyset$ for every $p \geq1$.
\end{example}

The recent developments obtained in \cite{BBDP} allow us to prove the
following characterization, which, besides its own interest, will be helpful
several times later.

\begin{theorem}
\label{caracterizacao} Let $P \in\mathcal{P}(^{n}E;F)$ and $a \in E$. Then $a
\in S_{p;q}(P)$ if and only if $0 \in S_{p;q}(\hat d^{k}P(a))$ for every $k =
1, \ldots, n$.
\end{theorem}

\begin{proof} Assume that $0 \in S_{p;q}(\hat d^kP(a))$ for every $k = 1, \ldots, n$. Given $(x_j)_{j=1}^\infty
\in \ell_q^w(E)$, $(\hat d^kP(a)(x_j))_{j=1}^\infty \in \ell_p(F)$ for $k = 1, \ldots, n$. For every $j$,
\begin{eqnarray}
P(a + x_j) - P(a) = \sum_{k=1}^{n}\binom{n}{k}\check P(a^{n-k},x_j^k) =  \sum_{k=1}^{n}\frac{1}{k!}\hat d^k
P(a)(x_j),\nonumber
\end{eqnarray}
so $a \in S_{p;q}(P)$ because $\ell_p(F)$ is a linear space. Conversely, let $a \in S_{p;q}(P)$ and $k \in
\{1, \ldots, n\}$. By \cite[Proposition 3.5]{BBDP} we know that $\check P$ is $(p;q)$-summing at $(a, \ldots,
a)$ in the sense of \cite[Section 2]{BBDP} or \cite[Definition 9.1]{BBJP}. So it follows from
\cite[Proposition 3.1]{BBDP} that the $k$-linear mapping
\[(x_1, \ldots, x_k) \in E^k \longrightarrow \check P(a^{n-k},x_1, \ldots, x_k)\in F\]
is $(p;q)$-summing at the origin. Calling on \cite[Proposition 3.5]{BBDP} once more it follows that the
polynomial generated by this $k$-linear mapping is $(p;q)$-summing at the origin. But this polynomial is a
multiple of $\hat d^kP(a)$, hence $\hat d^kP(a)$ is $(p;q)$-summing at the origin, that is $0 \in
S_{p;q}(\hat d^kP(a))$.
\end{proof}

First we apply this characterization to prove a substantial improvement of
Proposition \ref{Proposition 2.3}:

\begin{proposition}
\label{Proposition 3.3} Let $E$ be a Banach space and $n \geq2$. Then
$S_{1}(P) = E$ for every $P \in\mathcal{P}(^{n}E)$.
\end{proposition}

\begin{proof} Using multilinear mappings the result follows from an easy combination of \cite[Lemma 1]{irish} and \cite[Proposition 3.5]{BBDP}.
We prefer to provide a direct reasoning: let $P \in {\cal P}(^nE)$ and $a \in E$. For every $k = 1, \ldots,
n$, $\hat d^kP(a)$ is a scalar-valued $k$-homogeneous polynomial on $E$, so it is 1-summing at the origin by
Proposition \ref{Proposition 2.3}. The result follows from Theorem \ref{caracterizacao}.
\end{proof}

\begin{remark}
\rm Example \ref{Example 3.2} shows that Proposition \ref{Proposition 3.3} is no longer true for
vector-valued polynomials. Actually, we know much more: for every infinite-dimensional Banach space $E$,
every $p\geq1$, every $n\geq2$ and every $a\in E$, \cite[Theorem 3.7]{BBDP} assures that there is a
polynomial $P\in\mathcal{P}(^{n}E;E)$ such that $a\notin S_{p}(P)$.
\end{remark}

Following the same line of thought of Proposition \ref{Proposition 3.3} we obtain:

\begin{proposition}
Let $E$ be a Banach space with the Orlicz property. Then $S_{2,1}(P) = E$ for
every $n $, every $F$ and every $P \in\mathcal{P}(^{n}E;F)$.
\end{proposition}

\begin{proof} Given $(x_j)_{j=1}^\infty \in \ell_1^w(E)$, $(x_j)_{j=1}^\infty \in \ell_2(E)$ because $E$ has
the Orlicz property. Then,
\[\sum_{j=1}^\infty \|P(x_j)\|^2 \leq \left(\sum_{j=1}^\infty \|P(x_j)\|^{\frac{2}{n}}\right)^n
\leq \|P\|^2\left(\sum_{j=1}^\infty \|x_j\|^2\right)^n < + \infty.\]
This shows that $0 \in S_{2,1}(P)$ for every homogeneous polynomial $P$ on $E$. The result follows from
Theorem \ref{caracterizacao}.
\end{proof}

We finish this section showing that, for every $n$, the summing set of an
$n$-homogeneous polynomial may be $\{0\}$, thus may fail to be 1-lineable.

\begin{example}
\label{2-homogeneo}\rm Consider the 2-homogeneous polynomial
\[
P \colon L_{2}([0, 1];\mathbb{K}) \longrightarrow L_{1}([0, 1];\mathbb{K}%
)~;~P(f) = f^{2}.
\]
Let us see that $S_{1}(P) = \{0\}$. $0 \in S_{1}(P)$ by Proposition
\ref{orlicz} because $L_{2}([0, 1];\mathbb{K})$ has the Orlicz property. Let
$0 \neq f \in L_{2}([0, 1];\mathbb{K})$. Choose a sequence $(\alpha_{j})_{j
=1}^{\infty}$ in $\ell_{2}$ but not in $\ell_{1}$ and an orthonormal sequence
$(h_{j})_{j=1}^{\infty}$ in $L_{2}([0, 1];\mathbb{K})$ such that, for every $j
\in\mathbb{N}$, $|h_{j}(x)| = 1$ almost everywhere, Lebesgue measure (for
example, the Rademacher functions). Now we consider the sequence $(\alpha
_{j}h_{j})_{j=1}^{\infty}$. For every $g \in L_{2}([0, 1];\mathbb{K})$, by
Bessel's inequality we have
\begin{align}
\sum_{j = 1}^{\infty}|\langle g , \alpha_{j}h_{j} \rangle| \!\!  &  = \!\!
\sum_{j = 1}^{\infty}|\alpha_{j}| |\langle g , h_{j} \rangle| \leq\left(
\sum_{j = 1}^{\infty}|\alpha_{j}|^{2} \right)  ^{\frac{1}{2}} \cdot\left(
\sum_{j = 1}^{\infty}|\langle g , h_{j} \rangle|^{2} \right)  ^{\frac{1}{2}%
}\nonumber\\
\!\!  &  \leq\!\! \|(\alpha_{j})_{j =1}^{\infty}\|_{\ell_{2}}\|g\|_{L_{2}%
}.\nonumber
\end{align}
This shows that $(\alpha_{j}h_{j})_{j=1}^{\infty}\in\ell_{1}^{w}(L_{2}([0,
1];\mathbb{K}))$. On the other hand, since $\hat dP(f)(g) = 2 \check P(f,g) =
2 fg$, we have
\[
\sum_{j =1}^{\infty}\|\hat d P(f)(\alpha_{j}h_{j}) \|_{L_{1}} = 2\sum_{j =
1}^{\infty} \int_{0}^{1} |f(t)||\alpha_{j}||h_{j}(t)|dt = 2\|f\|_{L_{1}}
\sum_{j = 1}^{\infty}|\alpha_{j}| = + \infty,
\]
showing that $\hat dP(f)$ is not 1-summing. $f \notin S_{1}(P)$ by Theorem \ref{caracterizacao}.
\end{example}

Now we handle the case $n \geq3$. For each $x \in[- \pi, \pi]$, we consider
the set
\[
J_{x} := \{t \in[- \pi, \pi] : x-t \in[- \pi, \pi]\} = \left\{
\begin{array}
[c]{cl}%
[x - \pi,\pi], & 0 \leq x \leq\pi,\\
\mathrm{[}-\pi, x+\pi], & -\pi\leq x \leq0.
\end{array}
\right.
\]
Given $f, g \in L_{2}([-\pi, \pi];\mathbb{K})$, the \textit{convolution} $f *
g$ is defined on $[-\pi,\pi]$ by
\[
f * g(x) = \int_{J_{x}}f(x-t)g(t)dt.
\]
By Young's inequality \cite[Proposition 8.7]{folland} we know that $f*g \in
L_{2}([-\pi, \pi];\mathbb{K})$ and $\|f*g\|_{L_{2}} \leq\|f\|_{L_{2}%
}\|g\|_{L_{2}}$. Given $f \in L_{2}([-\pi, \pi];\mathbb{K})$ and $n \geq2$,
since the convolution is associative ($(f*g)*h = f*(g*h)$ - see
\cite[Proposition 8.6(b)]{folland}), inductively we can define
\[
f * \overset{(n)}{\cdots} * f := f * (f * \overset{(n-1)}{\cdots} * f).
\]

\begin{proposition}
\label{convolucao} Let $n \geq2$. Consider the $(n+1)$-homogeneous polynomial
\[
P \colon L_{2}([-\pi, \pi];\mathbb{K}) \longrightarrow L_{1}([-\pi,
\pi];\mathbb{K})~;~P(f) = (f * \overset{(n)}{\cdots} * f)\cdot f.
\]
Then, $S_{1}(P) = \{0\}$.
\end{proposition}

\begin{proof} $0 \in S_1(P)$ by Proposition \ref{orlicz} because $L_2([-\pi,
\pi];\mathbb{K})$ has the Orlicz property. Let $0 \neq f \in L_2([-\pi, \pi];\mathbb{K})$. For every $g \in
L_2([-\pi, \pi];\mathbb{K})$,
\begin{eqnarray}
\hat d P(f)(g) = (n+1) \check P(f^n,g) = (f * \stackrel{(n)}{\cdots} * f)\cdot g + n (f *
\stackrel{(n-1)}{\cdots} * f*g)\cdot f.
\end{eqnarray}
Claim 1: The linear operator
\[u \colon L_2([-\pi, \pi];\mathbb{K}) \longrightarrow L_1([-\pi, \pi];\mathbb{K})~;~u(g) = (f *
\stackrel{(n-1)}{\cdots} * f*g)\cdot f\] is 1-summing.\\
Proof of Claim 1: For $x \in [- \pi, \pi]$ define
\[h_x \colon [- \pi, \pi] \longrightarrow \mathbb{K}~;~h_x(t) =  \left\{\begin{array}{cl}
\overline{(f * \stackrel{(n-1)}{\cdots} * f)(x-t)}, & t \in J_x, \\ 0, & t \notin J_x.
\end{array} \right.\]
We have $h_x \in L_2([-\pi, \pi];\mathbb{K})$ and $\|h_x\|_{L_2} \leq \|f\|^{n-1}_{L_2}$. Regarding $h_x$ as
a linear functional on $L_2([-\pi, \pi]; \mathbb{K})$, for every $g \in L_2([-\pi, \pi];\mathbb{K})$ we
obtain
\begin{eqnarray}
h_x(g) \!\!& = &\!\! \langle g, h_x \rangle = \int_{J_x}g(t) \overline{h_x(t)}dt = \int_{J_x}g(t)(f *
\stackrel{(n-1)}{\cdots} * f)(x-t) dt\nonumber\\
\!\!& = &\!\! (f * \stackrel{(n-1)}{\cdots} * f*g)(x).\nonumber
\end{eqnarray}
We conclude that $u$ is 1-summing observing that, for $g_1, \ldots, g_k \in L_2([-\pi, \pi];\mathbb{K})$,
\begin{eqnarray}
\sum_{j=1}^k \|u(g_j)\|  \!\!& = &\!\! \sum_{j=1}^k \|(f * \stackrel{(n-1)}{\cdots} * f*g_j)\cdot
f\|_{L_1}\nonumber\\
\!\!& = &\!\! \sum_{j=1}^k \int_{-\pi}^\pi |(f * \stackrel{(n-1)}{\cdots} * f*g_j)(x)||f(x)| dx \nonumber\\
\!\!& = &\!\! \int_{-\pi}^\pi \sum_{j=1}^k |h_x(g_j)||f(x)| dx = \int_{-\pi}^\pi \|h_x\|_{L_2}\sum_{j=1}^k \left |\frac{h_x}{\|h_x\|_{L_2}}(g_j)\right ||f(x)|
dx\nonumber\\
\!\!& \leq &\!\! \|f\|^{n-1}_{L_2} \|f\|_{L_1}\|(g_j)_{j=1}^k\|_{w,1}.\nonumber
\end{eqnarray}
\medskip
\noindent Claim 2: The linear operator
\[v \colon L_2([-\pi, \pi];\mathbb{K}) \longrightarrow L_1([-\pi, \pi];\mathbb{K})~;~v(g) = (f * \stackrel{(n)}{\cdots} * f)\cdot g\]
fails to be 1-summing.\\
Proof of Claim 2: As we did in Example \ref{2-homogeneo}, let $(h_j)_{j=1}^\infty$ be an orthonormal sequence
in $L_2([-\pi, \pi];\mathbb{K})$ such that, for every $j \in \mathbb{N}$, $|h_j(x)| = \frac{1}{\sqrt{2\pi}}$
almost everywhere. Choosing a sequence $(\alpha_j)_{j = 1}^\infty$ in $\ell_2$ but not in $\ell_1$, the
argument we used in Example \ref{2-homogeneo} shows that $(\alpha_j h_j)_{j=1}^\infty \in \ell_1^w(L_2([-\pi,
\pi];\mathbb{K}))$. $v$ fails to be 1-summing because
\begin{eqnarray}
\sum_{j=1}^\infty \|v(\alpha_jh_j)\|  \!\!& = &\!\! \sum_{j=1}^\infty \|(f * \stackrel{(n)}{\cdots} * f)\cdot
\alpha_jh_j)\|_{L_1}\nonumber\\
\!\!& = &\!\! \sum_{j=1}^\infty \int_{-\pi}^\pi |(f * \stackrel{(n)}{\cdots} * f)(x)||\alpha_j||h_j(x)|dx\nonumber\\
\!\!& = &\!\! \frac{1}{\sqrt{2\pi}}\|(f * \stackrel{(n)}{\cdots} * f)\|_{L_1}\sum_{j=1}^\infty |\alpha_j| = + \infty.\nonumber
\end{eqnarray}
\medskip
\noindent Combining Claim 1, Claim 2 and (*) we conclude that $\hat d P(f)$ fails to be 1-summing. By Theorem
\ref{caracterizacao} it follows that $f \notin S_1(P)$.
\end{proof}

\section{Lineability and spaceability}

In this section we show that the case $n=3$ marks the loose of linearity:
while non-void summing sets of either $2$-homogeneous polynomials or
scalar-valued $3$-homogeneous polynomials are always linear subspaces, the
same is no longer true for vector-valued $3$-homogeneous polynomials. We also
show that for $n$-homogeneous polynomials, $n\geq4$, even in the scalar-valued
case the summing set may fail to be a linear subspace. We start by showing
that the summing set of a $2$-homogeneous polynomial is either empty or a
linear subspace:

\begin{theorem}
\label{Theorem 4.1} Let $E$ and $F$ be Banach spaces, $P \in\mathcal{P}%
(^{2}E;F)$ and $1 \leq q \leq p$. Then either $S_{p;q}(P) = \emptyset$ or
$S_{p;q}(P) = \{a\in E : \hat d P(a) \mathrm{~is~} (p;q)\mathrm{-summing}\}$.
So, $S_{p;q}(P)$ is either empty or a linear subspace of $E$. In particular,
if $P \in\mathcal{P}(^{2}E)$, then either $S_{p;q}(P) = \emptyset$ or
$S_{p;q}(P) = E$.
\end{theorem}

\begin{proof} Suppose $S_{p;q}(P) \neq \emptyset$. By Proposition \ref{Proposition 2.2} we have that
$P$ is $(p;q)$-summing at the origin. In this case, Theorem \ref{caracterizacao} yields that $a \in
S_{p;q}(P) \Longleftrightarrow \hat d P(a)$ is $(p;q)$-summing, which proves the first assertion. So,
$S_{p;q}(P) = (\hat d P)^{-1}(\Pi_{p;q}(E;F))$, which is a linear subspace of $E$ because $\hat d P$ is a
linear operator and $\Pi_{p;q}(E;F)$ is a linear subspace of ${\cal L}(E;F)$. For $P \in {\cal P}(^2E)$,
regardless of the vector $a \in E$, $\hat d P(a)$ is a linear functional, thus $(p;q)$-summing, so the last
assertion follows.
\end{proof}

Next example shows that, though always a linear subspace, the summing set of a
2-homogeneous polynomial may fail to be spaceable.

\begin{example}
\rm (A non-closed lineable non-spaceable summing set) Consider the 2-homogeneous polynomial
\[
P \colon\ell_{2} \longrightarrow\ell_{1} ~;~ P((\alpha_{j})_{j=1}^{\infty}) =
(\alpha_{j}^{2})_{j=1}^{\infty}.
\]
$P$ is 1-summing at the origin by Proposition \ref{orlicz} because $\ell_{2}$ has the Orlicz property. By
Theorem \ref{Theorem 4.1} it follows that $S_{1}(P) = \{a\in\ell_{2} : \hat d P(a) {\rm ~is~ 1-summing}\}$.
Given $a = (a_{k})_{k=1}^{\infty}\in\ell_{2}$, $\hat dP(a)$ is the linear operator
\[
(\alpha_{k})_{k=1}^{\infty}\in\ell_{2} \longrightarrow\hat d P(a)((\alpha
_{k})_{k=1}^{\infty}) = 2\check P((a_{k})_{k=1}^{\infty}), (\alpha_{k}%
)_{k=1}^{\infty})) = 2(a_{k} \alpha_{k})_{k=1}^{\infty}\in\ell_{1}.
\]
That is, $\frac{1}{2}\hat dP(a)$ is the diagonal operator by the vector $a$. By \cite[Theorem 9]{garling} it
follows that $\hat dP(a)$ is 1-summing if and only if $a \in\ell_{1}$, so $S_{1}(P) = \{a\in\ell_{2} : \hat d
P(a) {\rm~is~ 1-summing}\}= \ell_{1},$ which is a non-closed infinite-dimensional subspace of $\ell_{2}$
(obvious) that fails to be spaceable (it is well known that, as a subset of $\ell_{2}$, $\ell_{1}$ is not
spaceable).
\end{example}

Concerning summing sets, scalar-valued $3$-homogeneous polynomials behave like
$2$-homogeneous polynomials:

\begin{proposition}
Let $E$ be a Banach space and $P \in\mathcal{P}(^{3}E)$. Then, $S_{p;q}(P)$ is
either empty or a linear subspace of $E$ for every $1 \leq q \leq p $.
\end{proposition}

\begin{proof} Suppose $S_{p;q}(P) \neq \emptyset$ and let $a \in E$. By Proposition \ref{Proposition 2.2} we know that $P$ is $(p;q)$-summing at the origin. $\hat d P(a)$ is
$(p;q)$-summing because it is a linear functional. Calling on Theorem \ref{caracterizacao} it follows that $a
\in S_{p;q}(P)$ $\Longleftrightarrow$ $\hat d^2 P(a)$ is $(p;q)$-summing at the origin. We denote the space
of all scalar-valued 2-homogeneous polynomials on $E$ which are $(p;q)$-summing at the origin by ${\cal
P}_{as(p;q)}(^2E)$. Hence $S_{p;q}(P) = (\hat d^2 P)^{-1}({\cal P}_{as(p;q)}(^2E)$), which is a linear
subspace of $E$ because $\hat d^2P \colon E \longrightarrow {\cal P}(^2E)$ is a linear operator.
\end{proof}

Now we prove a multipurpose result:

\begin{theorem}
\label{Theorem 5.1} Let $E$ and $F$ be Banach spaces, $n \in\mathbb{N}$, $P
\in\mathcal{P}(^{n}E)$ and $g \colon E \longrightarrow F$ be a continuous
mapping. If $S_{p;q}(P) = E$, then $S_{p;q}(P\! \cdot g) = \mathrm{ker}P \cup
S_{p;q}(g)$.
\end{theorem}

\begin{proof} Given $(x_j)_{j=1}^\infty \in \ell_q^u(E)$, $g(x_j) \longrightarrow
g(0)$ because $x_j \longrightarrow 0$. It follows that $(g(x_j))_{j=1}^\infty$ is bounded, so there is $M
\geq 0$ such that $\|g(x_j)\| \leq M$ for every $j$. First let us show that $S_{p;q}(P\! \cdot g) \neq
\emptyset$. $(P(x_j))_{j=1}^\infty$ is absolutely $p$-summable because $S_{p;q}(P) = E$, so
\[\sum_{j=1}^\infty\|(P\! \cdot g)(x_j)\|^p = \sum_{j=1}^\infty|P(x_j)|^p\|g(x_j)\|^p \leq M^p\sum_{j=1}^\infty|P(x_j)|^p <
+ \infty,\] showing that $0 \in S_{p;q}(P\! \cdot g)$. Let $a \in E$. For every $j$,
\begin{eqnarray}
(P\! \cdot g)(a+x_j) - (P\! \cdot g)(a)  \!\!& = &\!\! P(a+x_j)g(a+x_j) - P(a)g(a)\nonumber\\
\!\!& = &\!\! P(a+x_j)(g(a+x_j) - g(a)) +\nonumber\\
\!\!& ~ &\!\!  + g(a)(P(a+x_j) - P(a)).\nonumber
\end{eqnarray}
We know that $((P(a+x_j) - P(a))_{j=1}^\infty$ is absolutely $p$-summable because $S_{p;q}(P) = E$, so $a \in
S_{p;q}(P\! \cdot g)$ if and only if $(P(a+x_j)(g(a+x_j) - g(a)))_{j=1}^\infty$ is absolutely $p$-summable.
$a \in S_{p;q}(P)$ because $S_{p;q}(P) = E$ by assumption, so Theorem \ref{caracterizacao} yields that $\hat
d^k P(a)$ is $(p;q)$-summing at the origin for every $k = 1, \ldots, n$. Combining this with the fact that
the sequence $(g(a+x_j) - g(a))_{j=1}^\infty$ is bounded (because $a + x_j \longrightarrow a$ and $g$ is
continuous), from
\begin{eqnarray}
P(a+x_j)(g(a+x_j) - g(a)) \!\!& = &\!\! P(a)(g(a+x_j) - g(a)) + (g(a+x_j) \nonumber\\
\!\!& ~ &\!\! - g(a))\sum_{k=1}^{n}\binom{n}{k}\check P(a^{n-k},x^k).\nonumber\\
\!\!& = &\!\! P(a)(g(a+x_j) - g(a)) + (g(a+x_j) \nonumber\\
\!\!& ~ &\!\! - g(a))\sum_{k=1}^{n}\frac{1}{k!}\hat d^kP(a)(x_j) {\rm ~for~every~}j,\nonumber
\end{eqnarray}
it follows that $a \in S_{p;q}(P\! \cdot g)$ if and only if $(P(a)(g(a+x_j) - g(a)))_{j=1}^\infty$ is
absolutely $p$-summable. So, ${\rm ker}P \subseteq S_{p;q}(P\! \cdot g)$ and for $a \notin {\rm ker}P$ we
have $a \in S_{p;q}(P\! \cdot g)$ if and only if $g$ is $p$-summing at $a$, that is $a \in S_{p;q}(g)$.
\end{proof}

A formula for $C_p(P\cdot Q)$ was proved in \cite[Theorem 5]{Aron-Dimant} for scalar-valued homogeneous
polynomials $P$ and $Q$. Observe that the formula we obtained above for $S_{p;q}(P\! \cdot g)$ holds true for
arbitrary continuous mappings $g$.

\begin{corollary}
\label{corolario 1} Let $E$ and $F$ be Banach spaces, $n \in\mathbb{N}$, $P
\in\mathcal{P}(^{n}E)$ and $g \colon E \longrightarrow F$ be a continuous
mapping. Then $S_{1}(P\! \cdot g) = \mathrm{ker}P \cup S_{1}(g)$.
\end{corollary}

\begin{proof} By Proposition \ref{Proposition 3.3} we know that $S_1(P) = E$. Now the result follows from Theorem
\ref{Theorem 5.1}.
\end{proof}

\begin{corollary}
\label{corollary}~\newline{\rm(a) (Complex case)} Let $E$ and $F$ be complex Banach spaces with $E$
infinite-dimensional. Let $P$ and $g$ be as in
Theorem \ref{Theorem 5.1}. Then $S_{1}(P\! \cdot g)$ is spaceable.\newline%
{\rm(b) (Real case)} Let $E$ and $F$ be real Banach spaces with $E$ infinite-dimensional. At least one of the
following possibilities occur:\newline\indent {\rm (b1)} There exists $P \in\mathcal{P}(^{3}E;E)$ such that
$S_{1}(P) = \{0\}$; \newline\indent {\rm (b2)} For every $P$ and $g$ as in Theorem \ref{Theorem 5.1} with $n
= 2$, $S_{1}(P\! \cdot g)$ is spaceable.
\end{corollary}

\begin{proof} (a) By Theorem \ref{Theorem 5.1} we know that ${\rm ker}P \subseteq S_1(P\! \cdot g)$, and from \cite[Theorem
5]{pz} we know that there exists an infinite-dimensional subspace
$G \subseteq {\rm ker}P$. But ${\rm ker}P = P^{-1}(0)$ is closed,
so $\overline{G} \subseteq \overline{{\rm ker}P} = {\rm ker}P
\subseteq S_1(P\! \cdot g)$.
\medskip
\noindent (b) Suppose that $E$ admits a positive definite 2-homogeneous polynomial $Q \in {\cal P}(^2E)$.
Defining $P \in {\cal P}(^3E;E)$ by $P(x) = Q(x)x$, by Theorem \ref{Theorem 5.1} we get $S_1(P) = {\rm ker}Q
= \{0\}$, proving that (b1) occurs in this case. If $E$ does not admit a positive definite 2-homogeneous
polynomial, by \cite[Theorem 1]{Pos} we know that for every $P \in {\cal P}(^2E)$, ${\rm ker}P$ contains an
infinite-dimensional subspace of $E$. A repetition of the proof of (a) shows that (b2) occurs in this case.
\end{proof}

\begin{remarks}
\rm (a) It is not always true that $kerP \subseteq S_{1}(P)$. For instance, if $P$ is the polynomial of
Example \ref{Example 3.2}, then $S_{1}(P) = \emptyset$ whereas $kerP = \{0\}$.\newline(b) It is interesting
to mention that possibility (b1) above occurs if there is a continuous linear injection from $E$ into a
Hilbert space, and possibility (b2) occurs otherwise (see \cite[Proposition 2]{Pos}).\newline(c) Corollary
\ref{corollary}(a) can be used to obtain information about non-reducibility of polynomials: if the polynomial
$P \in\mathcal{P}(^{n}E;F)$ between complex Banach spaces (dim$E = + \infty$) is such that $S_{1}(P)$ is
non-spaceable, then $P$ is irreducible, that is: $P$ cannot be written as $P = P_{1}\cdot P_{2}$ with $1 \leq
k \leq n-1$, $P_{1} \in\mathcal{P}(^{k}E)$ and $P_{2} \in\mathcal{P}(^{n-k}E;F)$. For example, the
convolution polynomials from $L_{2}([-\pi,\pi],\mathbb{C})$ to $L_{1}([-\pi,\pi],\mathbb{C}))$ of Proposition
\ref{convolucao} are irreducible.
\end{remarks}

\begin{example}
\label{exemplo 5.6}\rm(The summing set of a 3-homogeneous polynomial may fail to be a linear subspace) Let
$E$ be an infinite-dimensional Banach space. Fix $\varphi_{1}, \varphi_{2} \in E^{\prime}$ with
$ker\varphi_{1} \not \subseteq ker\varphi_{2}$ and $ker\varphi_{2} \not \subseteq ker\varphi_{1}$. For
example, choose linearly independent vectors $a, b \in E$ and functionals $\varphi_{1}, \varphi_{2} \in
E^{\prime}$ such that $\varphi_{1}(a) = \varphi_{2}(b) = 0$ and $\varphi_{2}(a) = \varphi_{1}(b) = 1$.
Consider the polynomial
\[
P \colon E \longrightarrow E~:~P(x) = \varphi_{1}(x)\varphi_{2}(x)x~;~P
\in\mathcal{P}(^{3}E;E).
\]
$S_{1}(id_{E}) = \emptyset$ because the identity operator on an
infinite-dimensional Banach space is never absolutely summing. From Corollary
\ref{corolario 1} we obtain
\[
S_{1}(P) = ker(\varphi_{1} \cdot\varphi_{2}) = ker\varphi_{1} \cup
ker\varphi_{2}.
\]
So, both $a$ and $b$ belong to $S_{1}(P)$. Assume for a while that $(a+b) \in S_{1}(P)$. So, $(a+b) \in
ker\varphi_{1}$ or $(a+b) \in ker\varphi_{2}$; and in this case we would have $b = (a+b)-a \in
ker\varphi_{1}$ or $a = (a+b)-b \in ker\varphi_{2}$ - a contradiction. We have just proved that $(a+b) \notin
S_{1}(P)$, therefore $S_{1}(P)$ is not a linear subspace of $E$. Note that $S_{1}(P)$ is spaceable, because
$ker\varphi_{1} \subseteq S_{1}(P)$.
\end{example}

\begin{proposition}
For every $n\geq4$ and every $2\leq q\leq p$ there exists a scalar-valued $n$-homogeneous polynomial $P$ so
that $S_{p;q}(P)\neq \emptyset$ and fails to be a linear subspace.
\end{proposition}

\begin{proof} Given $n \geq 4$ and $2 \leq q \leq p$, let $E$ be a Banach space which admit a polynomial $Q \in {\cal P}(^{n-2}E)$
such that $S_{p;q}(Q) = \emptyset$ (for example, the $(n-2)$-homogeneous polynomial on $\ell_2$ defined in
Example \ref{Example 3.1}). Let $\varphi_1, \varphi_2 \in E'$ be as in Example \ref{exemplo 5.6}, that is,
${\rm ker}\varphi_1 \not\subseteq {\rm ker}\varphi_2$ and ${\rm ker}\varphi_2 \not\subseteq {\rm
ker}\varphi_1$. Define $P := \varphi_1 \cdot \varphi_2 \cdot Q \in {\cal P}(^n E)$. Example \ref{tipo finito}
gives $S_{p;q}(\varphi_1 \cdot \varphi_2) = E$, so by Theorem \ref{Theorem 5.1} we have
\[S_{p;q}(P) = {\rm ker}(\varphi_1 \cdot \varphi_2) \cup S_{p;q}(Q) = {\rm ker}(\varphi_1 \cdot
\varphi_2)= {\rm ker}\varphi_1 \cup {\rm ker}\varphi_2,\] which fails to be a linear subspace of $E$.
\end{proof}

In all our examples and results thus far, the summing set of a homogeneous
polynomial is either void, $\{0\}$ or lineable. So, a final question concerns
the existence of non-trivial ($\neq\emptyset,\neq\{0\}$) non-lineable summing
sets. We shall solve this problem by proving that, given $n\geq2$, every
finite-dimensional subspace of $L_{2}$ is the summing set of a certain
$n$-homogeneous polynomial.\newline\indent Let $G$ be a complemented subspace
of a Banach space $E$. It is plain that the projection from $E$ onto $G$ is
$p$-summing if and only if $G$ is finite-dimensional. Nevertheless, for $q<p$,
the projection onto an infinite-dimensional complemented subspace may be
$(p;q)$-summing.

\begin{proposition}
\label{final} Let $G$ be a complemented subspace of $E$ such that the
projection from $E$ onto $G$ is $(p;q)$-summing. If there is a polynomial $P
\in\mathcal{P}(^{n}E;F)$ such that $S_{p;q}(P) = \{0\}$, $n \geq2$, then there
exists a polynomial $Q \in\mathcal{P}(^{n}E;F)$ such that $S_{p;q}(Q) = G$.
\end{proposition}

\begin{proof} Let
$H$ be the topological complement of $G$, that is $E = G \oplus H$. By $\pi_H, \pi_G \colon E \longrightarrow
E$ we denote the projections onto $H$ and $G$, respectively. Define $Q := P \circ \pi_H \in {\cal P}(^n E;
F).$ Let $a \in G$. Given $(x_j)_{j=1}^\infty \in \ell_q^w(E)$, $(\pi_H(x_j))_{j=1}^\infty \in \ell_q^w(E)$
because $\pi_H$ is a bounded linear operator. Hence $(Q(x_j))_{j=1}^\infty = (P(\pi_H(x_j)))_{j=1}^\infty \in
\ell_p(F)$ because $0 \in S_{p;q}(P)$. Since $\pi_H(a) = 0$, $(Q(a+ x_j)- Q(a))_{j=1}^\infty =
(Q(x_j))_{j=1}^\infty \in \ell_p(F)$, showing that $a \in S_{p;q}(Q)$. We proved that $G \subseteq
S_{p;q}(Q)$. Now we consider $a \notin G$. In this case $\pi_H(a) \notin S_{p;q}(P)$ because $\pi_H(a) \neq
0$. So we can find a sequence $(x_j)_{j=1}^\infty \in \ell_q^w(E)$ such that $(P(\pi_H(a)+ x_j)-
P(\pi_H(a)))_{j=1}^\infty \notin \ell_p(F)$. Since
\[P(\pi_H(a) + x_j) - P(\pi_H(a)) = \sum_{k=1}^{n}\binom{n}{k}\check P(\pi_H(a)^{n-k},x_j^k),\]
there is $k \in \{1, \ldots, n\}$ such that $(\check P(\pi_H(a)^{n-k},x_j^k))_{j=1}^\infty \notin \ell_p(F)$.
For every $x \in E$, $x = \pi_H(x) + \pi_G(x)$, so
\begin{eqnarray}
\check P(\pi_H(a)^{n-k},x_j^k) \!\!&=&\!\! \check P(\pi_H(a)^{n-k},(\pi_H(x_j) + \pi_G(x_j))^k)\nonumber\\
\!\!&=&\!\! \sum_{i=0}^{k}\binom{k}{i}\check P(\pi_H(a)^{n-k},\pi_H(x_j)^{k-i}, \pi_G(x_j)^i).\nonumber
\end{eqnarray}
Hence $(\check P(\pi_H(a)^{n-k},\pi_H(x_j)^{k-i}, \pi_G(x_j)^i))_{j=1}^\infty \notin \ell_p(F)$ for some $i
\in \{0, \ldots, k\}$. Assume, for a while, that $i \neq 0$. $\pi_G$ is $(p;q)$-summing by assumption, so
$\sum_{j=1}^\infty \|\pi_G(x_j)\|^p < + \infty$. Let $K$ be such that $\|x_j\| \leq K$ for every $j$. We have
\[\sum_{j=1}^\infty \|\check P(\pi_H(a)^{n-k},\pi_H(x_j)^{k-i},\pi_G(x_j)^i)\|^p \hspace*{15em}\]
\vspace*{-0.9em}
\[= \sum_{j=1}^\infty \|\check P(\pi_H(a)^{n-k},\pi_H(x_j)^{k-i},\pi_G(x_j)^{i-1},\pi_G(x_j) )\|^p \hspace*{2.5em}\]
\vspace*{-0.6em}
\[\hspace*{4.1em} \leq \|\check P\|^p
\|\pi_H\|^{(n-i)p}\|\pi_G\|^{(i-1)p}\|a\|^{(n-k)p}K^{(k-1)p}\sum_{j=1}^\infty\|\pi_G(x_j)\|^p < + \infty,\]
showing that $(\check P(\pi_H(a)^{n-k},\pi_H(x_j)^{k-i}, \pi_G(x_j)^i))_{j=1}^\infty \in \ell_p(F)$. It
follows that $i = 0$, that is $(\check P(\pi_H(a)^{n-k},\pi_H(x_j)^{k})_{j=1}^\infty \notin \ell_p(F)$. But
\[\check
P(\pi_H(a)^{n-k},\pi_H(x_j)^k) = (P \circ \pi_H)^\vee(a^{n-k},x_j^k) = \frac{(n-k)!}{n!}
\hat d^k Q(a)(x_j), \]
therefore $0 \notin S_{p;q}(\hat d^k Q(a))$. Now $a \notin S_{p;q}(Q)$ by Theorem \ref{caracterizacao}.
\end{proof}

\begin{corollary}
For every positive integer $n \geq2$ and every finite-dimensional subspace $G$
of $L_{2}([-\pi, \pi];\mathbb{K})$, there exists an $n$-homogeneous polynomial
$Q$ from $L_{2}([-\pi, \pi];\mathbb{K})$ to $L_{1}([-\pi, \pi];\mathbb{K})$
such that $S_{1}(Q) = G$.
\end{corollary}

\begin{proof} The projection from $L_2([-\pi, \pi];\mathbb{K})$ onto $G$ is 1-summing because it is a finite rank operator.
By Example \ref{2-homogeneo} and Proposition \ref{convolucao} we can consider a polynomial $P \in {\cal P}(^n
L_2([-\pi, \pi];\mathbb{K}); L_1([-\pi, \pi];\mathbb{K}))$ such that $S_1(P) = \{0\}$. The result follows
from Proposition \ref{final}
\end{proof}

\vspace*{1em} \noindent[Geraldo Botelho] Faculdade de Matem\'atica,
Universidade Federal de Uberl\^andia, 38.400-902 - Uberl\^andia, Brazil,
e-mail: botelho@ufu.br.

\medskip

\noindent[M\'ario C. Matos] Departamento de Matem\'atica, IMECC-UNICAMP, Caixa
Pos-tal 6065, 13.081-970 - Campinas, Brazil, e-mail: matos@ime.unicamp.br.

\medskip

\noindent\lbrack Daniel Pellegrino] Departamento de Matem\'{a}tica,
Universidade Federal da Para\'{\i}ba, 58.051-900 - Jo\~{a}o Pessoa, Brazil,
e-mail: dmpellegrino@gmail.com.

\end{document}